\documentclass{amsart}

\usepackage{amsmath,amsthm,amssymb,relsize}
\usepackage[margin=30mm]{geometry}
\usepackage{mathrsfs,mathtools,makecell}

\usepackage{booktabs, multirow, adjustbox, array}

\newcolumntype{R}[2]{%
    >{\adjustbox{angle=#1,lap=\width-(#2)}\bgroup}%
    l%
    <{\egroup}%
}

\usepackage{caption}
\usepackage{enumerate}
\usepackage{tikz}
\usetikzlibrary{cd, calc,  positioning, fit,  arrows, decorations.pathreplacing, decorations.markings, shapes.geometric, backgrounds,  bending
}
\usepackage{tikzsymbols}
\PassOptionsToPackage{dvipsnames}{xcolor}
\definecolor{darkblue}{rgb}{0,0,0.6}

\usepackage[breaklinks, pdftex, ocgcolorlinks,colorlinks=true, citecolor=darkblue, filecolor=darkblue, linkcolor=darkblue, urlcolor=darkblue]{hyperref}
\usepackage[capitalize,noabbrev]{cleveref}

\makeatletter
\newtheorem*{rep@theorem}{\rep@title}
\newcommand{\newreptheorem}[2]{%
\newenvironment{rep#1}[1]{%
 \def\rep@title{#2 \ref{##1}}%
 \begin{rep@theorem}}%
 {\end{rep@theorem}}}
\makeatother

\newtheorem{proposition}{Proposition}[section]
\newtheorem{theorem}[proposition]{Theorem}
\newtheorem{corollary}[proposition]{Corollary}
\newtheorem{lemma}[proposition]{Lemma}

\theoremstyle{definition}
\newtheorem{definition}[proposition]{Definition}

\newtheorem{example}[proposition]{Example}
\newtheorem{conjecture}[proposition]{Conjecture}

\theoremstyle{remark}

\newtheorem*{remark*}{Remark}

\newreptheorem{theorem}{Theorem}
\newreptheorem{lemma}{Lemma}
\newreptheorem{proposition}{Proposition}
\newreptheorem{corollary}{Corollary}
\newreptheorem{question}{Question}
\numberwithin{equation}{section}

\newcommand{\R}{{\ensuremath{\mathbb{R}}}}
\newcommand{\N}{{\ensuremath{\mathbb{N}}}}
\newcommand{\Z}{{\ensuremath{\mathbb{Z}}}}

\newcommand{\lk}{\operatorname{lk}}

\begin{document}

\title{A partial resolution of Hedden's conjecture on satellite homomorphisms}
%Other options:
%Patterns with even winding numbers not divisible by 8 do not induce concordance homomorphisms
\author{Randall Johanningsmeier}
\address{Department of Mathematics \& Statistics,  Swarthmore College, Swarthmore, PA USA}
\email{rjohann1@swarthmore.edu}

\author{Hillary  Kim}
\address{Department of Mathematics \& Statistics,  Swarthmore College, Swarthmore, PA USA}
\email{hkim11@swarthmore.edu}

\author{Allison N.\ Miller}
\address{Department of Mathematics \& Statistics,  Swarthmore College, Swarthmore, PA USA}
\email{amille11@swarthmore.edu}

\begin{abstract}
A pattern knot in a solid torus defines a self-map of the smooth knot concordance group. We prove that if the winding number of a pattern is even but not divisible by 8, then the corresponding  map is not a homomorphism, thus partially establishing a conjecture of Hedden.
\end{abstract}

\maketitle

\section{Introduction}

\begin{figure}[h!]
\[\begin{array}{ccc}
 \quad \begin{array}{c} \includegraphics[height=2.4cm]{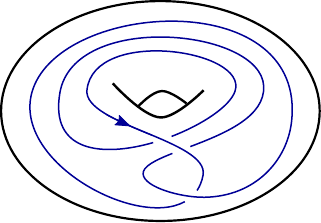} \end{array} \quad
& \quad \begin{array}{c} \includegraphics[height=2.7cm]{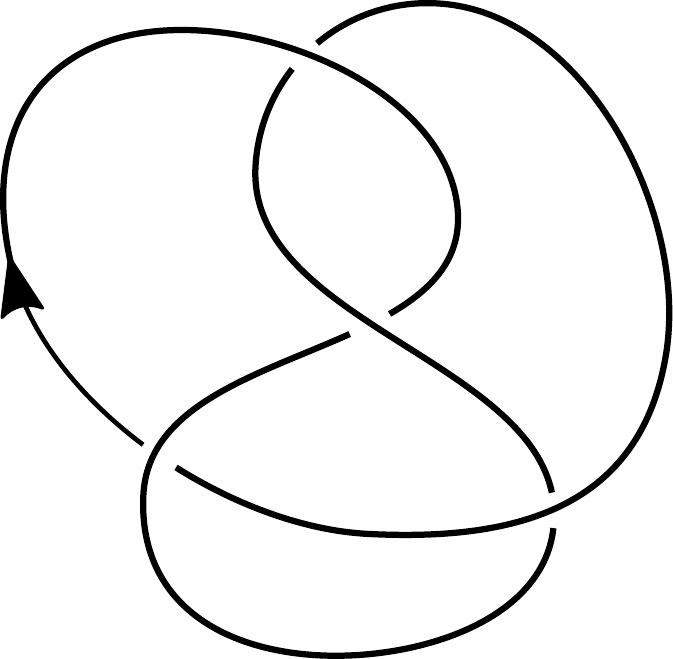} \end{array} \quad & \quad
\begin{array}{c} \includegraphics[height=2.7cm]{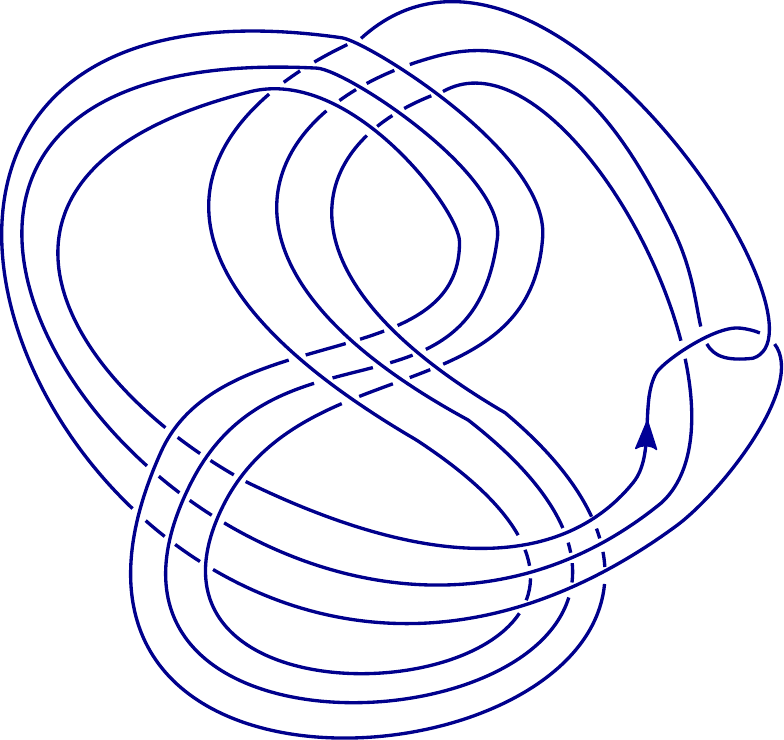} \end{array}
\end{array}\]
\caption{\small A pattern  $P$ (left), a knot $K$ (center), and the satellite knot $P(K)$ (right).} \label{fig:satellite}
\end{figure}

 The satellite construction 
 %illustrated in Figure~\ref{fig:satellite} 
 plays an important role in low-dimensional topology in general and knot concordance in particular. 
 %perhaps most notably in Thurston's classification of all knots as either embedding in a standard torus,  admitting a complete hyperbolic metric on their exterior, or resulting from a nontrivial satellite construction. 
A knot in a solid torus, or \emph{pattern}, induces a well-defined function on the set $\mathcal{C}$ of smooth concordance classes of knots, via the satellite operation illustrated in Figure~\ref{fig:satellite}. While these functions have been well-studied (see for example \cite{CochranHarveyGeometry, CochranDavisRay, CHLFractal,  DavisRay16, HeddenJPC, LevineMazur, KnotTraces, Ray15}), many questions remain open. 
In particular, $\mathcal{C}$ has the structure of an abelian group, with operation induced by connected sum.  While both the satellite operation and connected sum are geometrically defined operations,  they do not interact well: $P(K_1 \#K_2)$ is isotopic to $P(K_1) \# P(K_2)$ essentially only if $P$ is isotopic to either a core or an unknot in the solid torus. 

Our main result is the following, which roughly states that for patterns of certain (algebraic) winding number this behavior must persist even modulo concordance. \\

\begin{theorem}~\label{theorem:maintheorem}
Let $P$ be a pattern whose winding number  is even but not divisible by 8. Then $P$ does not induce a homomorphism of the smooth concordance group. 
\end{theorem}

This is progress towards establishing the following conjecture of Hedden. \\

\begin{conjecture} [\cite{BIRS16,MPIM16}]~\label{conj:hedden}
 Let $P$ be a pattern that induces a homomorphism of the smooth concordance group. 
 Then $P$ induces the zero map,  the identity map, or reversal. 
 \end{conjecture} 
It is straightforward to show that the zero map can only be induced by a winding number $0$ pattern, the identity map by a winding number $1$ pattern, and reversal by a winding number $-1$ pattern. In particular, Hedden's conjecture would imply that patterns whose winding number is greater than one in absolute value cannot induce homomorphisms. 
% We will henceforth refer to a pattern $P$ such that the induced map on smooth concordance is zero,  the identity, or reversal as \emph{standard}. 
Previous work in this area has obstructed specific examples like the Whitehead pattern~\cite{Gompf86} and $(n,1)$-cables~\cite{HeddenCabling} from inducing homomorphisms, as well as establishing new obstructions coming from Casson-Gordon signatures~\cite{Homomorphism} and the $d$-invariants of Heegaard Floer homology~\cite{LinkingNumberObstructions}. However, Theorem~\ref{theorem:maintheorem} is the first result to obstruct all patterns of a given winding number from inducing homomorphisms. Our proof relies on the following obstruction, due to recent work of Lidman, Miller, and Pinz{\'o}n-Caicedo. \\

 \begin{theorem}[\cite{LinkingNumberObstructions}]~\label{theorem:AJT}
 Let $P \subset S^1 \times D^2$ be a pattern of winding number $w$, where $w$ is divisible by some prime power $m$. 
 Let $\eta_m$ denote a preferred lift of $\eta=\{*\} \times \partial D^2$ to the $m$-fold branched cover $\Sigma_m(P(U))$, and let $t \colon \Sigma_m(P(U))) \to \Sigma_m(P(U))$ be a generator for the group of covering transformations. 
 
Suppose the following two conditions hold: 
\begin{enumerate}
    \item $\eta_m$ represents an odd order element of $H_1(\Sigma_m(P(U)))$,  
    \item  for $k=1, \dots, m-1$, the rational linking numbers $lk_{\Sigma_m(P(U))}(\eta_m ,t^k \eta_m)$ are not all equal to 0, while being either all non-positive or all non-negative. 
\end{enumerate}
Then there is a knot $K$ such that $P(-K) \# P(K)$ is not smoothly slice, and hence $P$ does not induce a homomorphism of the smooth concordance group.
 \end{theorem}

Note that $lk_{\Sigma_m(P(U))}(\eta_m ,t^k \eta_m)= lk_{\Sigma_m(P(U))}(\eta_m ,t^{m-k} \eta_m)$ for all $1 \leq k \leq m-1$, so to verify condition (2) of Theorem~\ref{theorem:AJT} we need consider at most $\lfloor \frac{m}{2} \rfloor$ linking numbers. 
Also, we will only use Theorem~\ref{theorem:AJT} when $m$ is a power of 2, in which case $|H_1(\Sigma_m(P(U)))|$ must be odd (see e.g.~\cite{SomeAspectsGordon}) and so condition (1) of Theorem~\ref{theorem:AJT} will be automatically satisfied. \\
 
To prove Theorem~\ref{theorem:maintheorem},
we will split into cases according to whether the winding number $w$ of $P$ satisfies $w \equiv 2$ mod 4 (Corollary~\ref{cor:n=2}) or $w \equiv 4$ mod 8 (Corollary~\ref{cor:n=4}). 
When $w \equiv 2$ mod 4, we argue that $lk_{\Sigma_2(P(U))}(\eta_2, t \eta_2)\neq 0$  in order to apply Theorem~\ref{theorem:AJT} with $m=2$. 
When $w \equiv 4$ mod 8, we first observe that if  $\lk_{\Sigma_2(P(U))}(\eta_2, t \eta_2) \neq 0$, then we can apply Theorem~\ref{theorem:AJT} with $m=2$.  
We prove that if $\lk_{\Sigma_2(P(U))}(\eta_2, t \eta_2)=0$, then  $lk_{\Sigma_4(P(U))}(\eta_4, t\eta_4)= 0$ (Proposition~\ref{prop:2fold4fold}), before showing that $lk_{\Sigma_4(P(U))}(\eta_4, t^2 \eta_4) \neq 0$  in order to apply Theorem~\ref{theorem:AJT} with $m=4$. In all of these cases our computation of linking numbers comes from relating an arbitrary winding number $w$ pattern to the $(w,1)$-cable pattern via crossing changes (Proposition~\ref{prop:crossingchangesPtoCn1}), lifting the surgery curves realizing these crossing changes to the appropriate branched cover, and comparing linking numbers there. \\

Finally, note that for every odd integer $w$, there is an example of a pattern $P$ with winding number $w$ such that $P$ is isotopic to $-P$, and hence such that $P(-K)$ is isotopic to $-P(K)$ for all knots $K$~\cite{Homomorphism}. 
However, there are no known examples of a pattern $Q$ with even winding number such that $Q(-K)$ is always concordant to $-Q(K)$, unless $Q$ induces the 0-map. This leads to the following conjecture.\\

\begin{conjecture}~\label{conj:pseudohom}
   Let $P$ be a pattern with nonzero even winding number. Then there exists a knot $K$ such that
$P(-K)$ is not concordant to $-P(K)$. 
\end{conjecture}
Since Theorem~\ref{theorem:AJT} obstructs a pattern $P$ from having the property that $P(-K)$ is always concordant to $-P(K)$, our proof of Theorem~\ref{theorem:maintheorem} establishes Conjecture~\ref{conj:pseudohom} for  patterns whose winding numbers are not divisible by 8. 

\subsection*{Acknowledgements}
Much of the work of this project took place while the first two authors were summer research students supported by the Natural Sciences \& Engineering division of Swarthmore College (RJ) and the Panaphil Foundation's Frances Velay Women's Science Research Fellowship (HK). 
We thank Tye Lidman and Arunima Ray for helpful comments on early versions of this paper.

\section{Definitions and examples}

The rational linking number of a pair of disjoint oriented simple closed curves in a rational homology sphere is defined as follows, see \cite[Chapter 10, Section 77]{SeifertThrelfall} for further details. \\

\begin{definition}~\label{definition:linking number}
    Let $\gamma_1$ and $\gamma_2$ be disjoint oriented simple closed curves in a rational homology 3-sphere $Y$. Let $\ell \in \mathbb{N}$ be such that $\ell\gamma_2$ represents the trivial element of $H_1(Y; \Z)$, and let $F$ be a 2-chain in $Y$ with boundary $\partial F= \ell \gamma_2$. 
    Then the \emph{(rational) linking number of $\gamma_1$ and $\gamma_2$} is 
    \[ lk_Y(\gamma_1, \gamma_2)= \frac{1}{\ell}( \gamma_1 \cdot F) \in \mathbb{Q}.\]
\end{definition}
Note that here and throughout the paper, for a simple closed curve $\gamma$ and a transverse 2-chain $F$,  we use $\gamma \cdot F$ to denote the signed count of intersection points between $\gamma$ and $F$. 
Although it is not particularly obvious from the definition, the linking number is symmetric and depends only on the homology class of one curve in the complement of the other. 

There is a natural correspondence between patterns in the solid torus and certain ordered 2-component links in $S^3$.
Given a pattern $P$ in $S^1 \times D^2$, we let $\eta=\{*\} \times \partial D^2$.  By considering the trivial embedding of $S^1 \times D^2$ in $S^3$, we obtain an ordered 2-component link $P(U) \cup \eta$ in  $S^3$, such that the second component $\eta$ is unknotted in $S^3$. In the other direction, given an ordered 2-component link $L_1 \cup L_2 \subseteq S^3$ with $L_2$ unknotted, we obtain a pattern by considering $L_1 \subseteq (S^3 \smallsetminus \nu(L_2)^\circ)\cong S^1 \times D^2$. We will frequently move between $P$ and $P(U) \cup \eta$ without much discussion. \\

\begin{example}~\label{exl:61}
The left side of Figure~\ref{fig:petatoetap} illustrates the 2 component link $P(U) \cup \eta$ describing the $(6,1)$ cable pattern, $C_{6,1}$. 
\begin{figure}[h!]
\begin{centering}
\begin{tikzpicture}
\node[anchor=south west,inner sep=0] at (0,0)
{\includegraphics[height=3cm]{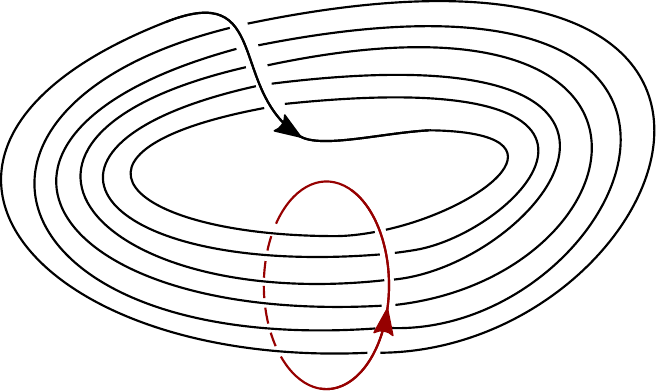}};
\node at (1,0.15){$P(U)$};
\node at (3, 0){$\eta$};
\end{tikzpicture}
\hspace{1cm}
\begin{tikzpicture}
\node[anchor=south west,inner sep=0] at (0,0)
{\includegraphics[height=3cm]{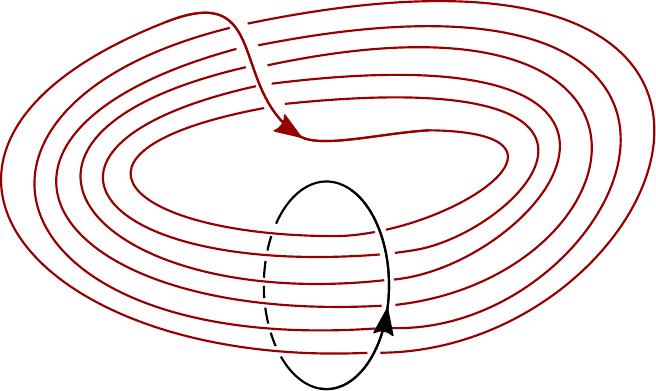}};
\node at (1.05,0.25){$\eta$};
\node at (3.2, -0.1){$P(U)$};
\end{tikzpicture}
\end{centering}
\caption{The link $P(U) \cup \eta$ defining the $C_{6,1}$ cable pattern is symmetric.}
\label{fig:petatoetap}
\end{figure}
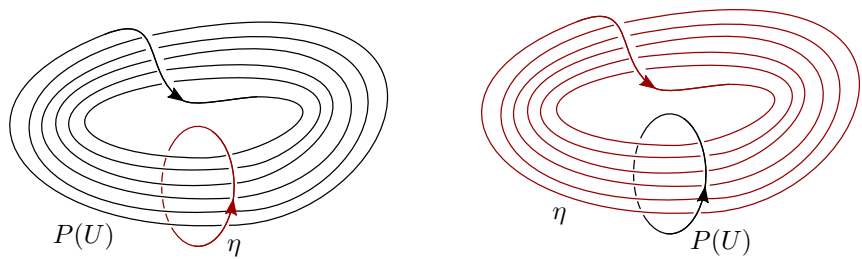
 On the right we have performed an isotopy so that $P(U)$ is the standard unknot. 
Now, observe that $\Sigma_2(P(U))= \Sigma_2(U)=S^3$, and the pre-image of $\eta$ in $\Sigma_2(P(U))$ is the $(6,2)$ torus link, whose components have linking number 3. 
\end{example}

This example generalizes to give the following. \\

\begin{lemma}~\label{lemma:cables}
Let $C_{n,1}(U) \cup \eta$ denote the 2-component link describing the $(n,1)$-cable pattern. Then
$lk_{\Sigma_m(C_{n,1}(U))}(\eta_m, t^j\eta_m)=\frac{n}{m}$ for all $m$ dividing $n$ and $1 \leq j \leq m-1$.
\end{lemma}

\begin{proof}
Suppose that $m$ divides $n$, and so $n=mk$ for some $k \in \N$. 
The link $C_{n,1}(U) \cup \eta$ is symmetric and the preimage of $\eta$ in $\Sigma_m(C_{n,1}(U))=S^3$ is the torus link $T(mk,m)$, as illustrated in Figure~\ref{fig:81} for $n=8$ and $m=4$.
\begin{figure}[h!]
    \begin{center}
     \includegraphics[height=2.5cm]{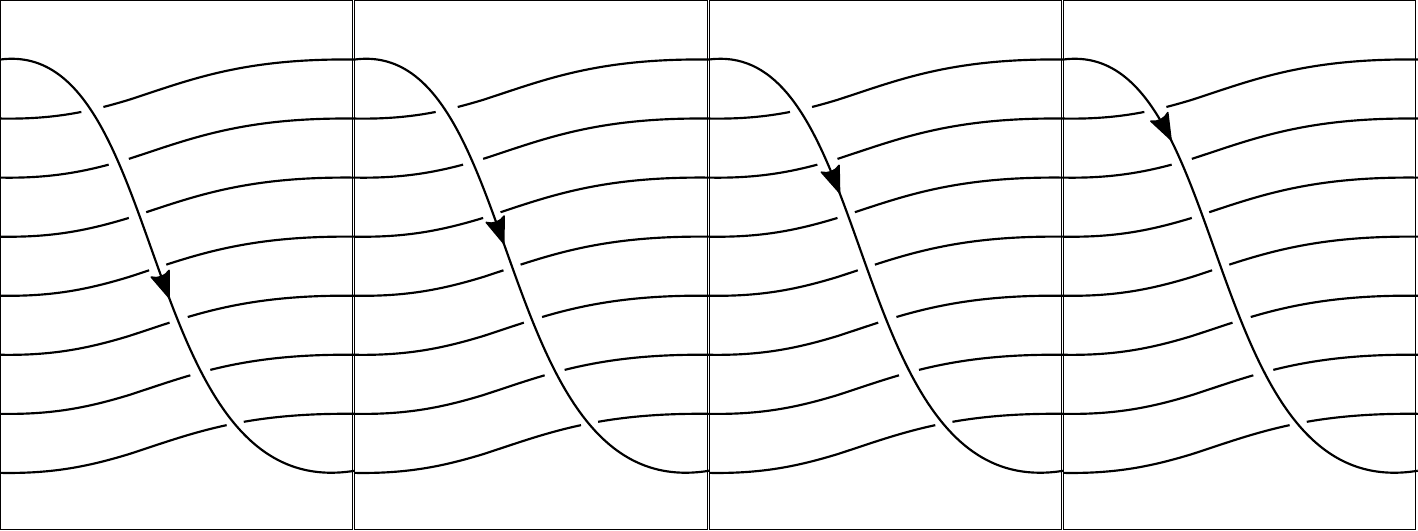}
    \end{center}
    \caption{The preimage of $\eta$ in $\Sigma_4(C_{8,1}(U))=S^3$, once the left and right sides of the diagram are identified  without twisting.}
    \label{fig:81}
\end{figure}
Since $T(mk,m)$ is isotopic to $T(m,mk)$, the $m$-component link obtained from the $m$-component unlink by inserting $k$ full twists between all components,  we see that the linking number between any two distinct lifts of $\eta$ is equal to $k= \frac{n}{m}$.  
\end{proof}

\begin{example}[A winding number 8 pattern that Theorem~\ref{theorem:AJT} does not obstruct from inducing a homomorphism.]
Consider the pattern depicted in Figure~\ref{fig:enter-label}.
    \begin{figure}[h!]
        \centering
        \begin{tikzpicture}
\node[anchor=south west,inner sep=0] at (0,0)
{\includegraphics[height=4cm]{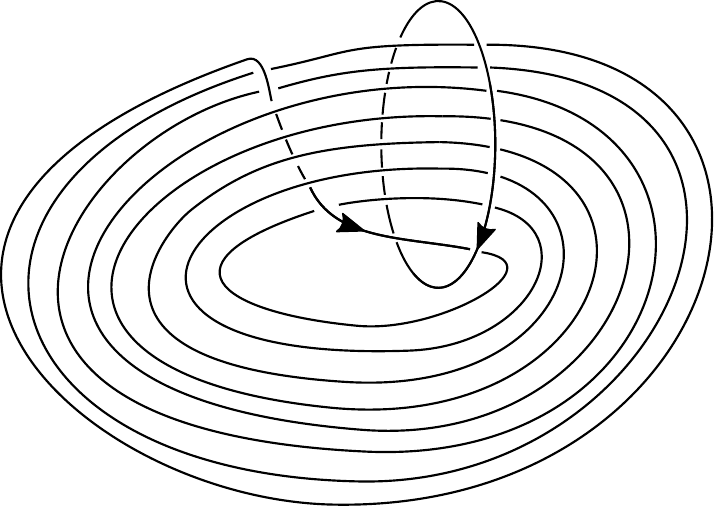}};
\node at (.65,0.2){$P(U)$};
\node at (3.9, 4){$\eta$};
\end{tikzpicture}
        \caption{A winding number 8 pattern that Theorem~\ref{theorem:AJT} does not obstruct from inducing a homomorphism.}
        \label{fig:enter-label}
    \end{figure}
The link $P(U) \cup \eta$ is symmetric, and since $P(U)=U$ and so $\Sigma_m(P(U))=S^3$ for all $m$, we can use a diagrammatic approach akin to that of Example~\ref{exl:61} and Lemma~\ref{lemma:cables} to compute the linking numbers corresponding to the $m$-fold cyclic branched covers for all prime powers $m$ dividing the winding number:  
\begin{itemize}
    \item $m=2$: $\lk_{S^3}(\eta_2, t\eta_2)=0$. 
    \item $m=4$:  $\lk_{S^3}(\eta_4, t\eta_4)=\lk_{S^3}(\eta_4, t^2\eta_4)=0$,
    \item $m=8$: $\lk_{S^3}(\eta_8, t\eta_8)=1$,  $\lk_{S^3}(\eta_8, t^2\eta_8)=0$, and $\lk_{S^3}(\eta_8, t^3\eta_8)=\lk_{S^3}(\eta_8, t^4 \eta_8)=-1$.
\end{itemize}
Since in each case the linking numbers are either identically zero or of mixed sign, Theorem~\ref{theorem:AJT} does not provide an obstruction.
\end{example}

Finally, we record the following for later use. 
\begin{proposition}\label{prop:crossingchangesPtoCn1}
Let $P$ be a pattern with winding number $n > 1$.  Then there exists a sequence of crossing changes which transforms $P$ into the cable pattern $C_{n, 1}$.
\end{proposition}
\begin{proof}
Note that a pattern with winding number $n$ represents the same homotopy class as $C_{n,1}$ in the solid torus.  As observed for example in~\cite[Section 4]{LevineHomotopyClassification}, a homotopy of curves can be realized by a composition of isotopies and crossing changes, thereby giving the desired result. The reader who finds this argument excessively terse or otherwise unsatisfying is referred to the appendix, where we give an alternate proof from a diagrammatic perspective. 
\end{proof}
In fact, an identical argument establishes Proposition~\ref{prop:crossingchangesPtoCn1} for any $n \in \Z$, so long as we interpret $C_{n,1}$ as a positively oriented core when $n=1$, an unknot in the solid torus when $n=0$, a negatively oriented core when $n=-1$, and as the reverse of $C_{|n|,1}$ when $n<-1$. 

% \subsection{Crossing changes to a cable}
% Our goal is now to establish the following proposition. 
% \begin{proposition}\label{prop:crossingchangesPtoCn1}
% Let $P$ be a pattern with winding number $n > 1$.  Then there exists a sequence of crossing changes which transforms $P$ into the cable $C_{n, 1}$.
% \end{proposition} We expect that this result is already known to the experts, see for example \cite[Lemma 3.5]{DNPR}, but since it is crucial for our proof of Theorem~\ref{theorem:maintheorem} we give a careful argument. 
% We also note that an analogous argument shows that every winding number 1 pattern can be transformed to a core of the solid torus via crossing changes and every winding number 0 pattern can be transformed to an unknot via crossing changes.

\section{Proof of main theorem}

The strategy for the proof of Theorem~\ref{theorem:maintheorem} is as follows. We will lift the  curves realizing the crossing changes converting the winding number $n$ pattern $P$ to the $(n,1)$ cable to the $m$-fold cyclic branched cover $\Sigma_m(C_{n,1}(U))=S^3$. This will give us a surgery description for $(\Sigma_m(P(U)), \cup_{i=0}^{m-1} t^i \eta_m^P)$ in terms of $(\Sigma_m(C_{n,1}(U)),\cup_{i=0}^{m-1} t^i \eta_m^C)$, where by a mild abuse of notation we use $t$ to refer to either of the covering transformation maps on $\Sigma_m(P(U))$ or on $\Sigma_m(C_{n,1}(U))$.
 The following theorem of Cha and Ko then allows us to relate the linking numbers of the components of $\cup_{i=0}^{m-1} t^i \eta_m^P$ in $\Sigma_m(P(U))$ to the linking numbers of the corresponding components of $\cup_{i=0}^{m-1} t^i \eta_m^C$ in $\Sigma_m(C_{n,1}(U))=S^3$, which we computed in Lemma~\ref{lemma:cables}.\\

\begin{theorem}[Theorem 3.1 of \cite{ChaKo02}]
~\label{theorem:chako}
Let $L=K_1 \cup \dots \cup K_\ell$ be an integrally framed link in $S^3$ such that the result of surgery on $S^3$ along $L$ is a rational homology 3-sphere $Y$. Let $A=(A_{i,j})$ be the linking-framing matrix of $L$. 
Then for any two disjoint 1-cycles $a$ and $b$ in $S^3 \smallsetminus L$, we have that 
\[lk_Y(a,b)= lk_{S^3}(a,b)-x^\intercal A^{-1}y,\]
where $x=(x_i)$ and $y=(y_i)$ are column vectors with $x_i= lk_{S^3}(a, K_i)$ and $y_i= lk_{S^3}(b,K_i)$ for all $i=1, \dots,\ell$. 
\end{theorem}

\subsection{The case of $n \equiv 2$ mod 4}

The following result relates the 2-fold cover linking numbers of a winding number $2k$ pattern to those of $C_{2k,1}$. \\

\begin{lemma}~\label{lemma:relatingtocable2}
Let $P=P(U) \cup \eta$ be a pattern with even winding number $n=2k$. 
Let $\eta^1$ and $\eta^2$ denote the lifts of $\eta$ to $\Sigma_2(P(U))$. 
There exists $a\in \mathbb{Z}$  such that 
\[ lk_{\Sigma_2(P(U))}(\eta_1, \eta_2)= k+ \frac{2a}{|H_1(\Sigma_2(P(U)))|}.\]
\end{lemma}

\begin{proof}
By Proposition~\ref{prop:crossingchangesPtoCn1}, $P$ can be transformed into the $(n,1)$ cable $C_{n,1}$ by a sequence of $m$ crossing changes for some $m \in \mathbb{N}$. We can realize the $i$th crossing change by doing $\epsilon_i$-framed surgery (for $\epsilon_i = \pm 1$) along a small curve $L_i$ linking $P(U)$ geometrically twice and algebraically zero times, while linking $\eta$  zero times. We refer to the resulting link  as $P_C \cup \eta_C \cup  \bigcup_{i=1}^m L_i$. Note that blowing down the $L_i$ curves transforms $P_C \cup \eta_C$ to $P \cup \eta$, while ignoring the $L_i$ curves we see $P_C \cup \eta_C= C_{n,1}$. 

Each framed curve $L_i$ lifts to a 2-component framed link $L_i^1 \cup L_i^2$ in $\Sigma_2(C_{n,1}(U))=S^3$, and the curve $\eta_C$ lifts to the 2-component link $\eta_C^1 \cup \eta_C^2$.  
Doing surgery along $L= \bigcup_{i=1}^m (L_i^1 \cup L_i^2)$ according to the induced framing converts $(\Sigma_2(C_{n,1}(U)), \eta_C^1 \cup \eta_C^2)$ to $(\Sigma_2(P(U)), \eta^1 \cup \eta^2). $

Let $A$ be the linking-framing matrix for $L$ with respect to the ordering $L_1^1, \dots, L_m^1, L_1^2, \dots L_m^2$.  Observe that since $lk(t \gamma_1, t \gamma_2)= lk(\gamma_1, \gamma_2)$ for any two curves $\gamma_1, \gamma_2$,
for any $1 \leq i, j \leq m$ we have
\begin{align*}
A_{i, j+m}&=lk(L_i^1, L_j^2)= lk(L_i^2, L_j^1)= A_{i+m,j} \\
A_{i+m, j+m}&= lk(L_i^2, L_j^2)= lk(L_i^1, L_j^1)= A_{i,j}. 
\end{align*}  
That is, there are $m \times m$ integer matrices $B$ and $C$ such that 
\[A= \begin{bmatrix}
    B & C \\ C & B
\end{bmatrix}, \]
i.e. $A$ is a $2 \times 2$ \emph{block circulant} matrix with size $m \times m$ blocks. 

Now define
\begin{align*}
    x&=(lk(\eta_C^1, L_1^1), \dots, lk(\eta_C^1, L_m^1), lk(\eta_C^1, L_1^2), \dots, lk(\eta_C^1, L_m^2))\\
    y&=(lk(\eta_C^2, L_1^1), \dots, lk(\eta_C^2, L_m^1), lk(\eta_C^2, L_1^2), \dots, lk(\eta_C^2, L_m^2)). 
\end{align*}
We have that $lk(\eta_C^2, L_i^2)= lk(t \eta_C^1, t L_i^1)=lk(\eta_C^1, L_i^1)$ and $lk(\eta_C^1, L_i^2)= lk(t \eta_C^2, t L_i^1)= lk(\eta_C^2, L_i^1)$ for all $i=1, \dots, m$. \\

\noindent \textbf{Claim: $lk(\eta_C^1, L_i^2)= - lk(\eta_C^1, L_i^1)$.}

Proof of claim: Let $F_i$ be a surface in $S^3$ with boundary $\partial F_i = L_i$. Observe that 
$0= lk(\eta_C, L_i)= F_i \cdot \eta_C.$
Let $\widetilde{F_i}$ denote a lift of $F_i$ to $\Sigma_2(C_{n,1}(U))=S^3$. 
Then $\widetilde{F_i}$ is a surface with boundary $\partial \widetilde{F_i}= L_i^1 \cup L_i^2$. Moreover, 
\begin{align*}
    \widetilde{F_i} \cdot (\eta_C^1 \cup \eta_C^2)= 2(F_i \cdot \eta_C)=0. 
\end{align*} 
Let $G_i$ be a surface in $\Sigma_2(C_{n,1}(U))=S^3$ with $\partial G_i= -L_i^1$. 
Then we can compute 
\begin{align*}
lk(\eta_C^1, L_i^2) + lk(\eta_C^2, L_i^2)
&= (\widetilde{F_i} \cup G_i) \cdot (\eta_C^1 \cup \eta_C^2)\\
&=  G_i \cdot (\eta_C^1 \cup \eta_C^2)
\\&= (G_i \cdot \eta_C^1)+ (G_i \cdot \eta_C^2)= - lk(\eta_C^1, L_i^1)- lk(\eta_C^2, L_i^1). 
\end{align*}
We can now apply the previous observations that $lk(\eta_C^2, L_i^1)= lk(\eta_C^1, L_i^2)$ and $lk(\eta_C^2, L_i^2)= lk(\eta_C^1, L_i^1)$ to finish the proof of our claim.  \checkmark\\

We therefore have that $x=(v, -v)$ and $y=(-v, v)$ for some $v \in \Z^m$. By Theorem~\ref{theorem:chako}, we know that 
\begin{align*}
lk_{\Sigma_2(P(U))}(\eta_1,\eta_2)&= lk_{S^3}(\eta_C^1, \eta_C^2) - x^TA^{-1}y= k-x^TA^{-1}y, 
\end{align*}
where we use Lemma~\ref{lemma:cables} to compute $lk(\eta_C^1, \eta_C^2)=k$.

It now only remains to show that $x^TA^{-1}y=\frac{2a}{|H_1(\Sigma_2(P(U)))|}$ for some $a \in \Z$. Since $A$ is a block circulant matrix, the main result of \cite{BlockCircular} implies that $A^{-1}$ is also block circulant, and hence can be written $A^{-1}= \begin{bmatrix} F & G \\ G & F\end{bmatrix}$ for some $F, G \in M_{m \times m} (\mathbb{Q})$.
We obtain
\begin{align*}
x^\intercal A^{-1} y= \begin{bmatrix} v^\intercal & -v^\intercal \end{bmatrix} \begin{bmatrix} F & G \\ G & F \end{bmatrix} \begin{bmatrix} -v \\v \end{bmatrix} 
= 2v^\intercal(G-F)v.
\end{align*}
The adjoint formula for a matrix inverse, applied to $A$, implies that the entries of $F$ and $G$ are all of the form $\frac{*}{\det(A)}$ for some $* \in \Z$, and  so we obtain our desired result, since $|\det(A)|=|H_1(\Sigma_2(P(U)))|$.
\end{proof}

The proof of Theorem~\ref{theorem:maintheorem} for $n \equiv 2 $ mod 4 now follows quickly.

\begin{corollary}\label{cor:n=2}[Theorem~\ref{theorem:maintheorem} for $n \equiv 2$ mod 4.]
    Let $P=P(U) \cup \eta$ be a pattern whose winding number is even but not divisible by 4.  Then $P$ does not induce a homomorphism of the concordance group. 
\end{corollary}
\begin{proof}
Write $n=2k$ for an odd integer $k$. 
By Lemma~\ref{lemma:relatingtocable2}, we know that there is $a \in \mathbb{Z}$ such that 
\[lk_{\Sigma_2(P(U))}(\eta_1, \eta_2)= k +\frac{2a}{|H_1(\Sigma_2(P(U)))|}= \frac{k |H_1(\Sigma_2(P(U)))|- 2a}{|H_1(\Sigma_2(P(U)))|} \neq 0,\]
since $k$ and $|H_1(\Sigma_2(P(U)))|$ are both odd.
Theorem~\ref{theorem:AJT} with $m=2$ then gives the desired result.
\end{proof}

\subsection*{The case of $n \equiv 4$ mod 8.}
Our proof strategy for $n \equiv 4$ mod 8 will involve  both the 2-fold and 4-fold branched cover linking numbers. The following result relates the two.  \\

\begin{proposition}~\label{prop:2fold4fold}
    Let $P= P(U) \cup \eta$ be a pattern with winding number that is divisible by 4. 
    Let $\eta_4^1$ denote a preferred lift of $\eta$ to $\Sigma_4(P(U))$, and let $\eta_4^{j+1}= t^j \eta_4^1$ for $j=1,2,3$. Define $\eta_2^1= \pi(\eta_4^1)$ and $\eta_2^2= \pi(\eta_4^2)$, where $\pi \colon \Sigma_4(P(U)) \to \Sigma_2(P(U))$ is the branched covering projection map. 
    Then \[lk_{\Sigma_2(P(U))}(\eta_2^1, \eta_2^2)= 2 lk_{\Sigma_4(P(U))}(\eta_4^1, \eta_4^2). \]
\end{proposition}

For convenience, in this proof we will abbreviate $lk_{\Sigma_2(P(U))}$ to $lk_2$ and $lk_{\Sigma_4(P(U))}$ to $lk_4$. 

\begin{proof}
Let $\ell=|H_1(\Sigma_4(P(U)))|$. Note that 
\[|H_1(\Sigma_4(P(U)))|= |\Delta_{P(U)}(-1)\Delta_{P(U)}(i)\Delta_{P(U)}(-i)|= |H_1(\Sigma_2(P(U)))| \cdot |\Delta_{P(U)}(i)\Delta_{P(U)}(-i)|,\]
where as usual $\Delta_{P(U)}(t)$ denotes the Alexander polynomial of $P(U)$. Therefore, every 1-cycle $a$ in $\Sigma_2(P(U))$ has the property that $\ell a$ bounds a 2-cycle.

Let $F$ be a 2-cycle in $\Sigma_2(P(U))$ such that $\partial F= \ell \eta_2^2$, so $lk_2(\eta_2^1, \eta_2^2)= \frac{1}{\ell}(\eta_2^1 \cdot F)$. 
Define $\widetilde{F}= \pi^{-1}(F)\subset \Sigma_4(P(U))$, and observe that $\partial \widetilde{F}= \ell(\eta_4^2 \cup \eta_4^4)$. Let $G$ be a 2-cycle in $\Sigma_4(P(U))$ such that $\partial G= -\ell \eta_4^4$. 
Observe that for $j \in \{1,3\}$ we have that
\begin{align*}
    lk_4(\eta_4^j, \eta_4^2)&= \frac{1}{\ell}(\eta_4^j \cdot (\widetilde{F} \cup G))= \frac{1}{\ell}(\eta_4^j \cdot \widetilde{F})+ \frac{1}{\ell}(\eta_4^j \cdot G)= 
    \frac{1}{\ell}(\eta_4^j \cdot \widetilde{F})- lk_4(\eta_4^j, \eta_4^4). 
\end{align*}
It follows that 
\begin{align*}
    4 lk_4(\eta_4^1, \eta_4^2)&=lk_4(\eta_4^1, \eta_4^2)+ lk_4(\eta_4^1, \eta_4^4)+ lk_4(\eta_4^3, \eta_4^2)+ lk_4(\eta_4^3, \eta_4^4)\\
    &= \frac{1}{\ell}((\eta_4^1 \cup \eta_4^3) \cdot \widetilde{F})\\
    &= \frac{1}{\ell}(\pi^{-1}(\eta_2^1) \cdot \pi^{-1}(F))\\
    &=\frac{2}{\ell} (\eta_2^1 \cdot F)\\
    &=2 lk_2(\eta_2^1, \eta_2^2), 
\end{align*}
where the second-to-last equality follows from the fact that $\pi \colon \Sigma_4(P(U)) \to \Sigma_2(P(U))$ is 2-to-1 everywhere besides the branch set $P(U)$.
\end{proof}

We also need the analogue of Lemma~\ref{lemma:relatingtocable2} for patterns with winding number divisible by 4. \\

\begin{lemma}~\label{lemma:relatingtocable4}Let $P= P(U) \cup \eta$ be a pattern of winding number $n=4k$.
Let $\eta^1, \eta^2, \eta^3, \eta^4$ denote the lifts of $\eta$ to $\Sigma_4(P(U))$, where for $i=1,2,3$, we obtain $\eta^{i+1}$  from $\eta^i$ by the action of the covering transformation.
Then there exists $a \in \mathbb{Z}$ such that 
\[ lk_{\Sigma_4(P(U))}(\eta^1, \eta^3)= k + \frac{2a}{|H_1(\Sigma_4(P(U)))|}.\]
\end{lemma}

\begin{proof}[Proof of Lemma~\ref{lemma:relatingtocable4}]
    Our strategy is extremely similar to that of the proof of Lemma~\ref{lemma:relatingtocable2}. Let $L_1 \cup \dots \cup L_m$ be a framed link of unknots,  surgery along which realizes the crossing changes of $P(U)$ that transform $P=P(U) \cup \eta$ to the $(n,1)$-cable pattern $C_{n,1}= P_C \cup \eta_C$. Note that  
    \[lk(L_i, P(U))=lk(L_i, \eta)=lk(L_i, L_j)=0
    \]for all $1 \leq i \neq j \leq m$

    For each $i=1, \dots, m$, pick a preferred lift $L_i^1$ of $L_i$ to $\Sigma_4(C_{n,1}(U))=S^3$, and let $L_i^2=tL_i^1$, $L_i^3= t^2 L_i^1$, and $L_i^4= t^3 L_i^1$. Similarly, let $\eta_C^1$ be a preferred lift of $\eta_C$ to $\Sigma_4(C_{n,1}(U))=S^3$, and let $\eta_C^2=t\eta_C^1$, $\eta_C^3= t^2 \eta_C^1$, and $\eta_C^4= t^3 \eta_C^1$.
    We have that $(\Sigma_4(P(U)), \cup_{i=1}^4 \eta^i)$ is obtained from $(\Sigma_4(C_{n,1}(U)), \cup_{i=1}^4 \eta_C^i)$ by performing appropriately framed surgery along 
    \[L= L_1^1 \cup \dots L_m^1 \cup L_1^2 \cup \dots L_m^2 \cup L_1^3 \cup \dots \cup L_m^3 \cup L_1^4 \cup \dots \cup L_m^4. 
    \]

    Now let $A$ be the linking-framing matrix of $L$  and let $x$ (respectively $y$) be the $4m$-component vector whose $i$th entry is the linking of $\eta_C^1$ (respectively $\eta_C^3$) with the $i$th component of $L$. Theorem~\ref{theorem:chako} tells us that 
    \[lk_{\Sigma_4(P(U))}(\eta^1, \eta^3)
    =lk(\eta_C^1, \eta_C^3)-x^TA^{-1}y= k-x^TA^{-1}y,
    \]
where for the last equality we use Lemma~\ref{lemma:cables}.

Observe that for any $1 \leq i, j \leq m$ and $1 \leq a, b \leq 4$ we have
\[lk(L_i^a, L_j^b)= lk(t^{a-1}L_i^1, t^{b-1}L_j^1)= lk(L_i^1, t^{b-a} L_j^1)= lk(L_i^1, L_j^{b-a+1}),
\]
where all exponents are taken modulo 4. 
It follows that $A$ is a block circulant matrix with blocks of size $m \times m$, and hence by the main result of \cite{BlockCircular} has a block circulant inverse 
\[A^{-1}=  \begin{bmatrix}
    Q&R&S&T\\T&Q&R&S \\ S&T&Q&R\\ R&S&T&Q
\end{bmatrix},\]
for some $Q,R,S,T \in M_{m \times m}(\mathbb{Q})$.
Note that since $A$ is a linking-framing matrix and hence symmetric, $A^{-1}$ is also symmetric, so  $Q^\intercal=Q$, $R^\intercal=T$, and $S^\intercal=S$.\\

\noindent \textbf{Claim:} $lk(\eta_C^1, L_i^4)=- lk(\eta_C^1, L_i^1)- lk(\eta_C^1, L_i^2)-lk(\eta_C^1, L_i^3)$.

Proof of claim:
Let $F_i$ be a surface in $S^3$ with boundary $\partial F=L_i$, and observe that $0=lk(\eta, L_i)= F \cdot \eta_C$. Lifting $F$ to $\Sigma_4(C_{n,1}(U))=S^3$, we obtain $\widetilde{F}$, a surface with boundary $L_i^1 \cup L_i^2 \cup L_i^3 \cup L_i^4$ and with the property that 
\[\widetilde{F} \cap (\eta_C^1 \cup \eta_C^3 \cup \eta_C^3 \cup \eta_C^4)= 4(F \cap \eta_C)=0.\]
Let $G_1$, $G_2$, and $G_3$ be surfaces in $\Sigma_4(C_{n,1}(U))=S^3$ with $\partial G_j=-L_i^j$, and observe that 
\begin{align*}
\sum_{k=1}^4 lk(\eta_C^k, L_i^4)&=
(\widetilde{F} \cup G_1 \cup G_2 \cup G_3) \cdot (\eta_C^1 \cup \eta_C^3 \cup \eta_C^3 \cup \eta_C^4) \\
&=( G_1 \cup G_2 \cup G_3) \cdot (\eta_C^1 \cup \eta_C^3 \cup \eta_C^3 \cup \eta_C^4)= \sum_{j=1}^3 \sum_{k=1}^4 -lk (\eta_C^k, L_i^j). 
\end{align*}
Rewriting, we obtain that 
\begin{align*}
    0&=(lk(\eta_C^1, L_i^1)+ lk(\eta_C^1, L_i^2)+ lk(\eta_C^1, L_i^3)+ lk(\eta_C^1, L_i^4)) \\
     &+(lk(\eta_C^2, L_i^2)+ lk(\eta_C^2, L_i^3)+ lk(\eta_C^2, L_i^4)+lk(\eta_C^2, L_i^1)) \\
     & +(lk(\eta_C^3, L_i^3)+ lk(\eta_C^3, L_i^4)+lk(\eta_C^3, L_i^1)+(lk(\eta_C^3, L_i^2)) \\
     &+ (lk(\eta_C^4, L_i^4)+lk(\eta_C^4, L_i^1)+lk(\eta_C^4, L_i^2)+lk(\eta_C^4, L_i^3)).
\end{align*}
Now observe that each 4-term parenthetical sum is equal, since $lk(\eta_C^j, L_i^k)$ depends only on $1 \leq i \leq m$ and the value of $k-j$ mod 4. So we obtain our desired claim. \checkmark \\

We therefore have that $x=(u,v,w,-u-v-w)$, where
$u_i= lk(\eta_C^1, L_i^1)$, $v_i=lk(\eta_C^1, L_i^2)$, and $w_i= lk(\eta_C^1, L_i^3)$ for $i=1, \dots, m$. 
Furthermore, since $x$ records the linking of $\eta_C^1$ with the components of $L$ and $y$ records the linking of $\eta_C^3=t^2 \eta_C^1$ with the components of $L$, we have that $y=(w, -u-v-w, u, v)$. 

We can now compute\begin{align*}
x^\intercal A^{-1} y &= 
\begin{bmatrix} u^\intercal & v^\intercal & w^\intercal & -u^\intercal-v^\intercal-w^\intercal \end{bmatrix}
\begin{bmatrix}
    Q&R&S&T\\T&Q&R&S \\ S&T&Q&R\\ R&S&T&Q
\end{bmatrix}
\begin{bmatrix}
    w\\ -u-v-w \\ u \\ v
\end{bmatrix}\\
&= u^\intercal Q w + v^\intercal Q (-u-v-w)+ w^\intercal Q u + (-u^\intercal-v^\intercal-w^\intercal) Q v \\
& \qquad + u^\intercal R (-u-v-w) + v^\intercal R u + w^\intercal R v +(-u^\intercal-v^\intercal-w^\intercal)Rw
\\ & \qquad \qquad +u^\intercal S u+v^\intercal S v+ w^\intercal Sw+ (-u^\intercal-v^\intercal-w^\intercal)S(-u-v-w)
\\ & \qquad  \qquad \qquad+ u^\intercal T v+v^\intercal T w+ w^\intercal T (-u-v-w)+
(-u^\intercal-v^\intercal-w^\intercal)Tu \\
&= 2(u^\intercal Q w - u^\intercal Q v - v^\intercal Q v - v^\intercal Q w-u^\intercal R u - u^\intercal R v - 2 u^\intercal R w - v^\intercal R w \\
& \qquad - w^ \intercal R w + v^\intercal R u + w^\intercal R v + u^\intercal S u + u^ \intercal S v+ u^\intercal S w+ v^\intercal Sv + v^\intercal S w+ w ^\intercal S w),
\end{align*}
where to obtain the final equality we repeatedly use that $Q^\intercal=Q$, $S^\intercal=S$, and $R^\intercal=T$. The entries of $u, v,$ and  $w$ are integers and the entries of $Q, R,$ and $ S$ are all of the form $\frac{*}{\det(A)}$ for $* \in \Z$, by the adjoint formula for a matrix inverse. Therefore, since $|\det(A)|=|H_1(\Sigma_4(P(U)))|$, we have our desired result. 
\end{proof}

Theorem~\ref{theorem:maintheorem} for $n \equiv 4$ mod 8 follows quickly from this result together with Proposition~\ref{prop:2fold4fold}. 
\begin{corollary}[Theorem~\ref{theorem:maintheorem} for $n \equiv 4$ mod 8.]
\label{cor:n=4}  
    Let $P$ be a pattern with winding number $n$ equivalent to 4 mod 8. Then $P$ does not induce a homomorphism of the concordance group. 
\end{corollary}

\begin{proof}
Write $n=4k$ for some odd $k$. 
If $\lk_{\Sigma_2(P(U))}(\eta_2, t \eta_2) \neq 0$, then by applying Theorem~\ref{theorem:AJT} with $m=2$ we are done. So assume that $\lk_{\Sigma_2(P(U))}(\eta_2, t \eta_2) = 0$. By Proposition~\ref{prop:2fold4fold}, this implies that $\lk_{\Sigma_4(P(U))}(\eta_4, t \eta_4) = 0$ as well. However, by Lemma~\ref{lemma:relatingtocable4}, we know that there is $a \in \Z$ such that 
\[\lk_{\Sigma_4(P(U))}(\eta_4, t^2 \eta_4)=k+ \frac{2a}{|H_1(\Sigma_4(P(U)))|} \neq 0,
\]
since $|H_1(\Sigma_4(P(U)))|$ is odd. Theorem~\ref{theorem:AJT} with $m=4$ then gives the desired result.
\end{proof}

\section*{Appendix}

In this section we give a diagrammatic proof that every winding number $n>1$ pattern can be converted to the $(n,1)$ cable by crossing changes.  We will parameterize the solid torus as
\[ S^1 \times D^2 = S^1 \times I_y \times I_z= \{(\theta, y,z): \theta \in S^1=([0,1]/\sim),\, y \in [0,1],\, z \in [0,1]\}.
\]
Given a pattern $P$ in the solid torus, a small isotopy ensures that $P$  intersects the disc $\theta=1$ transversally in some finite number of points, the quantity of which, counted without sign, we call the \emph{wrapping number} of this embedding of $P$.\footnote{This differs slightly from the usual definition of the wrapping number, which involves a minimum. Note that if we counted with sign we would obtain $n$, the winding number of $P$.}
By cutting open the solid torus along the $\theta=1$ disc, we obtain a tangle in $I_\theta \times I_y \times I_z \subseteq \R^3$. We call this a \emph{cube decomposition} of $P$, and each component of the cube decomposition is called a \emph{strand} of the pattern.

\begin{proof}[Proof of Proposition~\ref{prop:crossingchangesPtoCn1}]
By a small isotopy of $P$, we can assume that there is a unique point $p$ with minimal $y$-coordinate, as illustrated on the left of Figure~\ref{fig:crossingchanges}. 
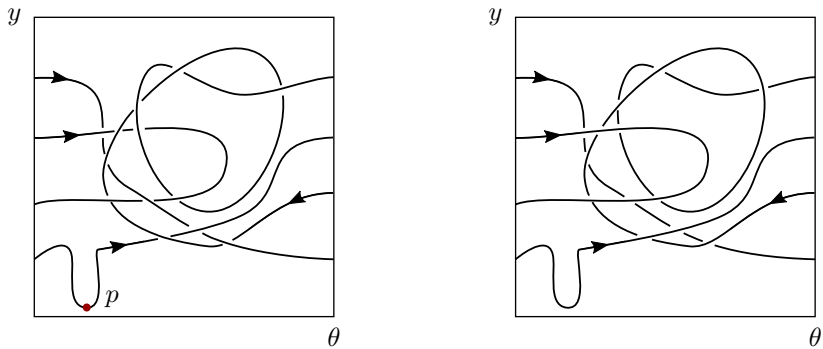
\begin{figure}[h!]
\begin{centering}
\begin{tikzpicture}
%\draw[step=1cm,color=gray] (0,0) grid (4,4);%Uncomment this to get some helpful grid lines
\node[anchor=south west,inner sep=0] at (0,0)
{\includegraphics[height=4cm]{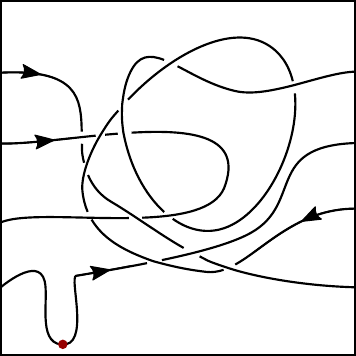}};
\node at (1.05,0.25){$p$};
\node at (-0.25,4){$y$};
\node at (4, -0.25){$\theta$};
\end{tikzpicture}
    \hspace{1.5cm}
    \begin{tikzpicture}
\node[anchor=south west,inner sep=0] at (0,0)
{\includegraphics[height=4cm]{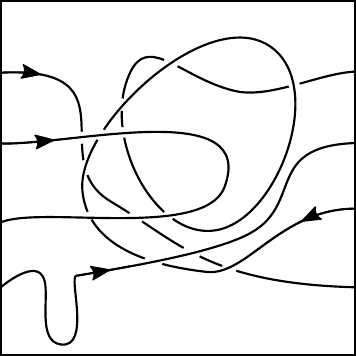}};
\node at (-0.25,4){$y$};
\node at (4, -0.25){$\theta$};
\end{tikzpicture}
    \end{centering}
    \caption{$P$ (left) is altered by crossing changes to give $Q$ (right).}
    \label{fig:crossingchanges}
\end{figure}
% \begin{figure}[h!]
% \begin{center}
% \includegraphics[height=8cm]{ycoordinateb.pdf}
% \end{center}
% \caption{$p$ is the only point on $P$ with $y$ coordinate $b$; every other point has $y$ coordinate greater than $b$.}
% \label{fig:ycoordinateb}
% \end{figure}
Now  consider the projection of $P$ onto the $\theta y$-plane. Beginning at $p$, move along $P$ according to its orientation, changing crossings as necessary so that every crossing is first reached as an overcrossing. We call this new pattern $Q$.
Observe that we can assign $z$-coordinates to $Q$ so that $p$ has $z$-coordinate $1-\varepsilon$ for some small $\epsilon>0$ and, as we move along $Q$, the $z$-coordinate is monotonically decreasing to $\varepsilon$ until just before we return to $p$, when it rapidly increases back to $1-\varepsilon$  along a small unknotted arc. 
We  say that embeddings such that the $z$-coordinate of the pattern is monotonically decreasing along the pattern except near one point, have the ``downhill property." We will now perform a series of isotopies while preserving the downhill property, in order to recognize $Q$ as $C_{n,1}$ or $C_{n,-1}$. 

First, suppose that the wrapping number of $Q$ is greater than its winding number. There must therefore be two strands in the cube decomposition that are \emph{returning}: that is, they enter and exit on the same side of the diagram.  Thus, we can select a returning strand that does not contain the point $p$, as illustrated on the left of Figure~\ref{fig:pullingaround}.
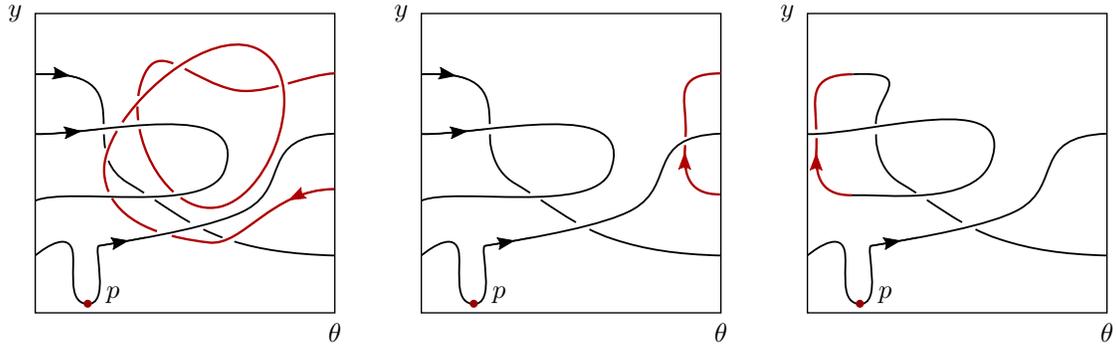
\begin{figure}[h!]
\begin{center}
\begin{tikzpicture}
\node[anchor=south west,inner sep=0] at (0,0)
{\includegraphics[height=4cm]{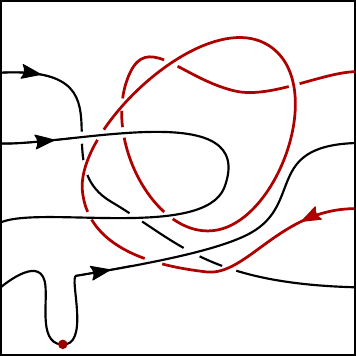}};
\node at (1.05,0.25){$p$};
\node at (-0.25,4){$y$};
\node at (4, -0.25){$\theta$};
\end{tikzpicture}
\quad
\begin{tikzpicture}
\node[anchor=south west,inner sep=0] at (0,0)
{\includegraphics[height=4cm]{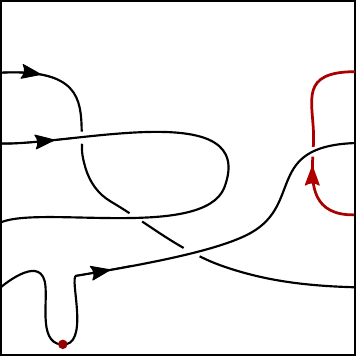}};
\node at (1.05,0.25){$p$};
\node at (-0.25,4){$y$};
\node at (4, -0.25){$\theta$};
\end{tikzpicture}
\quad
\begin{tikzpicture}
\node[anchor=south west,inner sep=0] at (0,0)
{\includegraphics[height=4cm]{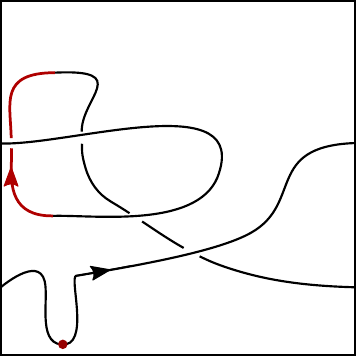}};
\node at (1.05,0.25){$p$};
\node at (-0.25,4){$y$};
\node at (4, -0.25){$\theta$};
\end{tikzpicture}
\end{center}
\caption{Reducing wrapping number by two by isotoping a returning strand.}
\label{fig:pullingaround}
\end{figure}
Such a strand has monotonically decreasing $z$-coordinates that are distinct from the $z$-coordinates of every other point in the diagram, besides the small unknotted arc near $p$. Therefore we can isotope this strand to be without crossings and lie in a small neighborhood of the right side of the diagram, as illustrated in the center of Figure~\ref{fig:pullingaround}. 
A small isotopy around the torus then gives us a diagram whose wrapping number is strictly two less than the the original wrapping number, as on the right of Figure~\ref{fig:pullingaround}. Throughout these isotopies, we both preserve the $z$-coordinate of every point, hence preserving the downhill property,  and ensure that $p$ remains the point with minimal $y$-coordinate. 

After repeating this procedure as necessary, we can therefore assume that our diagram of $Q$  has wrapping number equal to winding number. Note that every strand of the resulting diagram has monotonically decreasing $z$-coordinate, except in a neighborhood of $p$. We now isotope each strand not containing $p$ to the straight line segment connecting its endpoints, while isotoping the strand containing $p$ to the union of three line segments: a $z$-decreasing line segment ending at a point $p'$ just before $p$, a line segment from $p'$ to $p$ with minimal $y$-value and $z$-value from $z=\varepsilon$ to $z=1-\varepsilon$ (we will call this the \emph{ascending segment}, and another $z$-decreasing line segment starting at $p$. As before, these isotopies can be simultaneously performed while preserving the $z$-coordinate of each point except for those on the ascending segment and so that no points besides those on the ascending segment have minimal $y$-coordinate. An example of the result of this process when $n=4$ is shown on the left of  Figure~\ref{fig:swapprojection}, where we have labeled the $z$-coordinates of the end of each line segment.
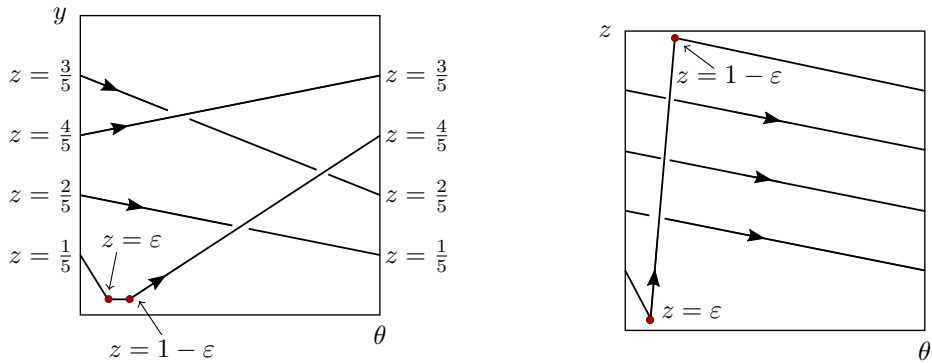
\begin{figure}[h!]
    \begin{center}
    \begin{tikzpicture}
\node[anchor=south west,inner sep=0] at (0,0)
{\includegraphics[height=4cm]{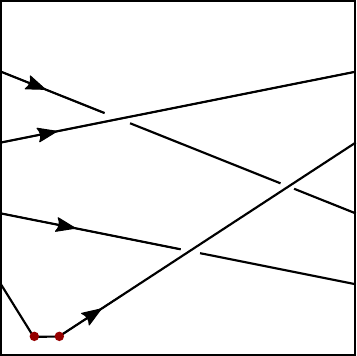}};
\node at (-0.25,4){$y$};
\node at (4, -0.25){$\theta$};
\node at (-0.5, 0.8){$z=\frac{1}{5}$};
\node at (-0.5, 1.6){$z=\frac{2}{5}$};
\node at (-0.5, 2.4){$z=\frac{4}{5}$};
\node at (-0.5, 3.2){$z=\frac{3}{5}$};
\node at (4.5, 0.8){$z=\frac{1}{5}$};
\node at (4.5, 1.6){$z=\frac{2}{5}$};
\node at (4.5, 2.4){$z=\frac{4}{5}$};
\node at (4.5, 3.2){$z=\frac{3}{5}$};
\node at (1.1, -0.45){$z=1-\varepsilon$};
\node at (.7, 1){$z=\varepsilon$};
\draw [-to] (1.1,-0.25)--(.8, .2);
\draw [-to] (.5,.85)--(.4,.35);
\end{tikzpicture}
\hspace{1.5cm}
\begin{tikzpicture}
\node[anchor=south west,inner sep=0] at (0,0)
{\includegraphics[height=4cm]{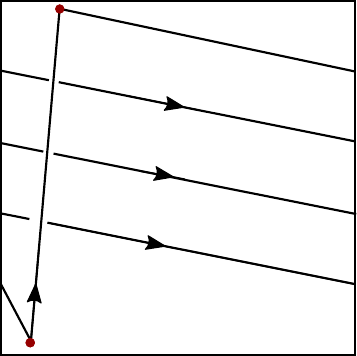}};
\node at (-0.25, 4){$z$};
%\node at (0.65,0.25){$p'$};
%\node at (0.9,3.6){$p$};
\node at (4, -0.25){$\theta$};
\node at (.9, 0.25){$z=\varepsilon$};
\node at (1.4,3.4){$z=1-\varepsilon$};
\draw [-to] (1.2,3.5)--(.8, 3.8);
\end{tikzpicture}    \end{center}
    \caption{An example of projection of $Q$ onto the $\theta y$- and $\theta z$-planes.}
    \label{fig:swapprojection}
\end{figure}

Now consider the projection of $Q$ onto the $\theta z$-plane, as illustrated on the right of Figure~\ref{fig:swapprojection}. Because this embedding of $Q$ has the downhill property, the only way for two points on $Q$ to have the same $z$ coordinate is if one of them lies on the ascending segment. Therefore, all crossings include the ascending segment. However, because $p$ has the highest $z$ coordinate and $p'$ has the lowest $z$ coordinate of all points on $Q$, and because all strands have been pulled taut while keeping $p$ and $p'$ fixed, by the intermediate value theorem  the ascending segment must cross each other strand exactly once in this projection. The points on the ascending segment have minimal $y$ coordinate $b$, so all of those crossings have the $p'$-to-$p$ line segment as the overstrand.\footnote{The reader might have expected `understrand', but we remind them of the orientation of $I_\theta \times I_y \times I_z \subseteq \R^3$.}

The diagram we have just described is a standard picture of either $C_{n, 1}$ or $C_{n,-1}$. If we have $C_{n, 1}$, we are done. If we have $C_{n, -1}$, we need only start the whole process again, still using the same point $p$ as a starting point, but moving against the orientation of $P$ as we change crossings so that every crossing is the opposite from in $Q$. 
\end{proof}
\bibliography{bib}
\bibliographystyle{plain}

\end{document}